\documentclass[12pt,leqno]{article}

\usepackage{amssymb, graphicx, amscd, amsmath, amssymb, amsthm, array, fancyvrb, relsize, paralist}

\newtheorem{thm}{Theorem}
\newtheorem{lem}[thm]{Lemma}

\numberwithin{equation}{section}

\def\F{\mathbb{F}}
\def\fx4{\frac{x}{4}}

\def\3F2{{}_3\hspace{-1pt}F_2}
\def\2F1{{}_2\hspace{-1pt}F_1}
\def\h2F1{{}_2\hspace{-1pt}\widehat{F}_1}
\def\Fq2{\F_{q^2}}
\def\e{\varepsilon}
\def\c4{\chi_4}
\def\oc4{\overline{\chi_4}}
\def\oA{\overline{A}}
\def\oB{\overline{B}}
\def\oC{\overline{C}}
\def\oD{\overline{D}}

\def\l({\left(}
\def\r){\right)}
\def\bar{\begin{array}{r|}}
\def\ear{\end{array}}

\title{A QUADRATIC HYPERGEOMETRIC $\2F1$ 
TRANSFORMATION OVER FINITE FIELDS}

\author{\\ \\ Ron Evans\\
Department of Mathematics\\
University of California at San Diego\\
La Jolla, CA  92093-0112\\
revans@ucsd.edu
\\ \\
and
\\ \\
John Greene\\
Department of Mathematics and Statistics\\
University of Minnesota--Duluth\\
Duluth, MN  55812\\
jgreene@d.umn.edu
}

\medskip

\date{May 18, 2016}

\begin{document}
\maketitle

\noindent 2010 \textit{Mathematics Subject Classification}.
11T24, 33C05.

\noindent \textit{Key words and phrases}.
Hypergeometric $\2F1$ functions over finite fields, Gauss sums, Jacobi sums,
pseudo hypergeometric functions, quadratic transformations, Hasse--Davenport
relation.

\newpage

\begin{abstract}
In 1984, the second author conjectured a quadratic transformation formula which
relates two hypergeometric $\2F1$ functions over a finite field $\F_q$.  
We prove this conjecture in Theorem 2.  The proof depends on
a new linear transformation formula for pseudo hypergeometric 
functions over $\F_q$. Theorem 2 is then applied to give an elegant new
transformation formula (Theorem 3) for $\2F1$ functions over finite fields.
\end{abstract}

\maketitle

\section{Introduction}
Let $\F_q$ be a field of $q$ elements, 
where $q$ is a power of an odd prime $p$. 
Throughout this paper, $A, B, C, D,  \chi,  \e, \phi$ 
denote complex multiplicative characters on $\F_q^*,\ $ 
extended to map 0 to 0. 
Here $\e$ and $\phi$ always denote the trivial and quadratic characters, 
respectively. For $y \in \F_q$, 
let $\zeta^y$ denote the additive character 
\begin{equation*}
\zeta^y := 
\exp \Bigg( \frac{2 \pi i}{p} \Big( y^p + y^{p^2} + \dots + y^q \Big) \Bigg).
\end{equation*}
Recall the definitions of the Gauss sum
\begin{equation*}
G(A) = \sum_{y \in \F_q} A(y) \zeta^y
\end{equation*}
and the Jacobi sum
\begin{equation*}
J(A, B) = \sum_{y \in \F_q} A(y) B(1-y).
\end{equation*}
These sums have the familiar properties 
\[
G(\e) = -1, \quad J(\e,\e) = q-2,
\]
and for nontrivial $A$, 
\[
G(A) G(\oA) = A(-1) q, \quad J(A, \oA) = -A(-1), 
\quad J(\e, A)=-1. 
\]
Gauss and Jacobi sums are related by \cite[p. 59]{BEW}
\begin{equation*}
J(A,B) = \frac{G(A) G(B)}{G(AB)}, \quad \text{if } AB \neq \e.
\end{equation*}
The Hasse--Davenport product relation \cite[p. 351]{BEW} yields
\begin{equation}
A(4) G(A) G(A \phi) = G(A^2) G(\phi).
\end{equation}

As in \cite[p. 82]{TAMS}, 
define the hypergeometric $\2F1$ function over $\F_q$ by
\begin{equation}
\2F1 \l( \bar A,B \\ C \ear \ x \r)
=\frac{\e (x)}{q}
\sum_{y \in \F_q} B(y) \oB C(y-1)
\oA(1-xy), \quad x \in \F_q.
\end{equation}
Define the binomial coefficient over $\F_q$
as in \cite[p. 80]{TAMS} by
\begin{equation*}
\begin{pmatrix}
A \\ B
\end{pmatrix}
= \frac{B(-1)}{q} J(A, \oB).
\end{equation*}
In \cite[(1.11)]{EG2}and \cite[(1.10)]{EG1}, 
we defined a pseudo hypergeometric function
$F^\ast (C,D; x)$ for $x \in \F_q$ by 
\begin{equation}
F^\ast(C,D; x) := \frac{q}{q-1} \sum_\chi
\begin{pmatrix} C \chi^2 \\ \chi \end{pmatrix}
\begin{pmatrix} C \chi \\ D \chi \end{pmatrix}
\chi \l( \fx4 \r) + CD(-1) \frac{\oC(x/4)}{q},
\end{equation}
where the sum is over all characters $\chi$ on $\F_q$.
The function $F^\ast (C,D; x)$ could be defined more simply
by using Gauss sums instead of Jacobi sums; see \cite[Lemma 2.1]{EG2}.
(This is in line with the observation
in \cite{DM} about an elegant way to define general 
hypergeometric functions over finite fields.)
However, we stick here with the formulation in (1.3) because it allows us
to more easily apply theorems in \cite{TAMS}.

We will need the following alternative formula for $F^\ast (C,D; x)$.
In \cite[p. 224]{EG2}, it is proved that if $C \neq D$ and 
$x \notin \{0,1\}$, then 
\begin{equation}
F^*(C,D;x) = \frac{C(2)}{q} \sum_t C \oD^2(1-t) \oC D(1-x-t^2).
\end{equation}
In fact, (1.4) holds even when $C=D$, but this fact will not be
used here, so we omit the proof.

\noindent 
{\em{Note}}:
On the fifth line from the bottom in \cite[p. 224]{EG2},
the misprint $\oA C^2$ should be corrected to $\oA C$.
Four lines before Theorem 1.1 in \cite{EG2}, replace the clause
``because ..." with ``because every element in $\F_q$
is a square in $\F_{q^2}$", and in the same sentence, replace
``(3.2)--(3.3)" with ``(2.16)--(2.17)".

\medskip

We are now prepared to state our results.

\begin{lem}
Suppose that $A$, $A^2\oB$, and $\phi A \oB$ are all nontrivial.
Let $y \in \F_q$ with $y \notin \{0,1\}$.  Then
\[
F^\ast(B, \oA B; y) =
\frac{\phi A B(-1) \oA^2 B(2) G(A^2 \oB)G(\phi \oA B)}{G(\phi) G(A)}
F^\ast(B, \phi A; 1 - y). 
\]
\end{lem}

Lemma 1 gives a linear transformation formula for
the pseudo hypergeometric function $F^\ast$.  It will be employed
to prove  Theorem 2.

\begin{thm}
Suppose that $A$, $A^2\oB$, and $\phi A \oB$ are all nontrivial.
Let $x \in \F_q$ with $x \neq -1$.  Then
\begin{equation} 
\begin{split}
& \2F1 \l( \bar A,B \\ A^2 \ear \ \frac{4x}{(1+x)^2} \r) = \\
& \frac{\oA(4) \phi B(-1)G(A^2 \oB)G(\phi \oA B)}{G(\phi) G(A)} B^2(1+x)\ 
\2F1 \l( \bar \phi \oA B,B \\ \phi A \ear \ x^2 \r).
\end{split}
\end{equation}
\end{thm}

Theorem 2 gives a finite field analogue
of an important $\2F1$ quadratic transformation of Gauss related to
elliptic integrals
\cite[p. 50]{Berndt}, \cite[(3.1.11)]{AAR}. 
A transformation equivalent to Theorem 2 is given in
\cite[Theorem 17]{JFLL}.

In 1984,
the second author \cite[(4.40)]{Thesis} proved Theorem 2 in the
special case that the character $B$ is even, and he conjectured that
Theorem 2 holds in general  \cite[p. 54]{Thesis}.  After proving
this conjecture, we will employ Theorem 2 to prove Theorem 3.

\noindent
{\em{Note}}:
On the second line of \cite[p. 54]{Thesis},
the misprint $2/(1+x)^2$ should be corrected to $2/(1+x^2)$.  Also, the
second equality in \cite[(4.40)]{Thesis} should be ignored, as it is incorrect.

\begin{thm}
Let $q \equiv 1 \pmod 4$, so that there exists a quartic character
$\c4$ on $\F_q$.  Let $z \in \F_q$ with $z \notin \{0,1,-1\}$.  Then
for any character $D$ on $\F_q$,
\begin{equation}
D^4(z-1) \2F1 \l( \bar D,D\chi_4 \\ \chi_4 \ear \ z^4 \r)
=\2F1 \l( \bar D,D^2\phi \\ D\phi \ear \ -\left(\frac{z+1}{z-1} \right)^2 \r).
\end{equation}
\end{thm}

Theorem 3, which motivated this paper, gives an elegant
transformation formula for $\2F1$ functions over finite fields.
The first author \cite{Evans} has applied Theorem 3
to evaluate a weighted sum of hypergeometric functions over $\F_q$.
The evaluation turns out to be elementarily equivalent to an identity that
Katz \cite{Katz} had proved using 
rigidity properties of Kloosterman sheaves.

Stanton \cite{Stant} has found the following analogue
of (1.6) over the complex numbers, valid for any non-negative integer $n$.
\begin{equation}
\begin{split}
(z-1)^{4n+2} &\2F1 \l( \bar -n-1/4,-2n-1 \\ -n+1/4 \ear
\ -\left(\frac{z+1}{z-1} \right)^2  \r)
\equiv \\
&(-2z)\frac{\Gamma(2n+3)\Gamma(3/4)}{\Gamma(n+2)\Gamma(n+3/4)}
\2F1 \l( \bar -n-1/4, -n \\ 5/4 \ear \ z^4  \r).
\end{split}
\end{equation}
Both sides of (1.7) are polynomials in $z$ of degree $4n+1$, and
the symbol $\equiv$ signifies that the two polynomials are identical.

\medskip

\noindent
{\em{Remark}}:
Transformation formulas for hypergeometric functions over $\F_q$
have numerous applications to number theory, algebraic geometry,
and modular forms; for some recent examples, see
\cite{BG}, \cite{JF}, \cite{JFDM}, \cite{Len}, \cite{DM}, \cite{DMMP}.
In \cite{JFLL}, a number of such transformation formulas
are proved and interpreted geometrically.

\section{Proof of Lemma 1}

In view of (1.4) and the hypothesis that $A$ and $\phi A \oB$ are
nontrivial, it suffices to prove that $\alpha = \beta$, where
\[
\alpha:= G(A^2 \oB)G(\phi \oA B)\sum_t \oA^2B(1-t) \phi A \oB (y-t^2)
\]
and 
\[
\beta = \phi B A(-1) A^2 \oB(2) G(\phi)G(A) \sum_t A^2 \oB(1-t) \oA(1-y-t^2).
\]
We have 
\[
\alpha = \sum_t \sum_w \sum_z \oA^2 B(1-t) A^2 \oB(w)
\phi A \oB(y-t^2) \phi \oA B(z) \zeta^{w+z}.
\]
Since $A^2\oB$ and $\phi \oA B$ are nontrivial, 
\[
\alpha = \sum_t \sum_{w \neq 0} \sum_{z \neq 0}
A^2 \oB(w) \phi  \oA B(z) \zeta^{w(1-t) +z(y-t^2)}.
\]
Replacing $w$ by $2wz$, we obtain
\begin{equation*}
\begin{split}
\alpha &= A^2 \oB(2) \sum_w A^2 \oB(w) \sum_z \phi A(z)
\zeta^{z(y-1+(w+1)^2)} \sum_t \zeta^{-z(t+w)^2} \\
&= A^2 \oB(2) \phi(-1) G(\phi)\sum_w A^2 \oB(w) 
\sum_z A(z) \zeta^{z(y-1+(w+1)^2)}.
\end{split}
\end{equation*}
Since $A$ is nontrivial, it follows that
\[
\alpha= A^2 \oB(2) \phi(-1) G(\phi)G(A) \sum_w A^2 \oB(w) \oA( y-1+(w+1)^2).
\]
Replacing $w$ by $w-1$, we see that $\alpha = \beta$,
which completes the proof of Lemma 1.

\section{Proof of Theorem 2}

Both sides of (1.5) vanish when $x=0$.  When $x=1$,
each $\2F1$ in (1.5) has the argument $1$, so that (1.5) can be 
directly verified using \cite[Theorem 4.9]{TAMS}, with the 
aid of the Hasse-Davenport relation (1.1).  Thus for the remainder of the
proof, assume that $x \notin \{0, -1, 1\}$.

Applying \cite[Theorem 4.4(i)]{TAMS} to the left side of (1.5), 
we see that (1.5)
is equivalent to
\begin{equation} 
\begin{split}
& \2F1 \l( \bar A,B \\ \oA B \ear \ \frac{(1-x)^2}{(1+x)^2} \r) = \\
& \frac{\oA(4) \phi A B(-1)G(A^2 \oB)G(\phi \oA B)}{G(\phi) G(A)} B^2(1+x)\
\2F1 \l( \bar  \phi \oA B,B \\ \phi A \ear \ x^2 \r).
\end{split}
\end{equation}

To prove (3.1), first suppose that $x=\pm i$, where $i \in \F_q$ is a primitive fourth root of unity. In this case $q \equiv 1 \pmod 4$ and
$\phi(-1)=1$. Each $\2F1$ in (3.1) has the argument
$-1$, so that (3.1) can be verified using \cite[(4.11)]{TAMS}, with the
aid of the Hasse-Davenport relation (1.1).  Thus for the remainder of the
proof, we assume that $x \notin \{0, -1, 1, i, -i\}$.

Apply \cite[Theorem 4.16]{TAMS} 
to see that the $\2F1$ on the left side of (3.1) equals
\[
\oB\left(\frac{2(x^2+1)}{(x+1)^2}\right) 
F^\ast\left(B, \oA B; \left(\frac{x^2-1}{x^2+1}\right)^2\right),
\]
by (1.3).  Similarly, the $\2F1$ on the right side of (3.1) equals
\[
\oB(1+x^2) 
F^\ast\left(B, \phi A; \frac{4x^2}{(1+x^2)^2}\right).
\]
Thus (3.1) is equivalent to
\begin{equation*}
\begin{split}
& F^\ast\left(B, \oA B; \left(\frac{x^2-1}{x^2+1}\right)^2\right) =\\
& =\frac{\phi A B(-1) \oA^2 B(2) G(A^2 \oB)G(\phi \oA B)}{G(\phi) G(A)}
F^\ast\left(B, \phi A; \frac{4x^2}{(1+x^2)^2}\right),
\end{split}
\end{equation*}
which immediately follows from Lemma 1 with 
\[
y=\left(\frac{x^2-1}{x^2+1}\right)^2.
\]
This completes the proof of Theorem 2.

\section{Proof of Theorem 3}

Note that since $q \equiv 1 \pmod 4$, we have $\phi(-1)=1$,
and there exists a primitive fourth root of unity 
$i \in \F_q$. The proof may be facilitated by the observation
that $-4 = (1+i)^4$, so that $\c4(-4)=1$.

If $D$ is either trivial or quartic, then the $\2F1$ functions
in (1.6) are degenerate, so that (1.6) follows directly
from \cite[Corollary 3.16]{TAMS}.  Thus assume that $D$ is
neither quartic nor trivial.

We will apply three transformations to convert the right side of (1.6)
to the left side.
First apply the transformation in \cite[Theorem 4.20]{TAMS} 
with $A=D$ and $B=D\c4$
to express the $\2F1$ on the right side of (1.6) in terms of
\begin{equation}
\2F1 \l( \bar D\oc4,D\c4 \\ D\phi \ear \ -\frac{(z^2-1)^2}{4z^2} \r).
\end{equation}
We will next utilize the transformation
\begin{equation}
\2F1 \l( \bar A, B \\ C \ear \ x  \r) =
ABC(-1)\oB(x)\2F1 \l( \bar B\oC, B \\ B\oA \ear \ \frac{1}{x}  \r),
\end{equation}
which follows by replacing $y$ by $y/x$ in (1.2).
Apply the transformation (4.2) with $A=D\oc4$, $B=D\c4$, and $C=D\phi$
to express the $\2F1$ in (4.1) in terms of
\begin{equation}
\2F1 \l( \bar \oc4,D\c4 \\ \phi \ear \ -\frac{4z^2}{(z^2-1)^2}  \r).
\end{equation}
Finally apply the transformation in Theorem 2 of this paper 
with $A=\oc4$, $B=D\c4$, and $x = -z^2$
to express the $\2F1$ in (4.3) in terms of the $\2F1$ on the 
left side of (1.6).  After some simplication,
these three successive transformations yield the desired result (1.6).

\end{document}